\documentclass[12pt,leqno,twoside]{article}
\usepackage{amssymb}
\usepackage{amsmath}
\usepackage{amsthm}
\usepackage{a4,indentfirst,latexsym}
\usepackage{graphics}
\usepackage{mathrsfs}
\usepackage{cite,enumitem,graphicx}
\usepackage[colorlinks=true,urlcolor=black,
citecolor=black,linkcolor=black,linktocpage,pdfpagelabels,
bookmarksnumbered,bookmarksopen]{hyperref}
\usepackage{color}

\newtheorem{thm}{Theorem}[section]
\newtheorem{prop}[thm]{Proposition}

\newtheorem{cor}[thm]{Corollary}
\newtheorem{rem}[thm]{Remark}

\newenvironment{altproof}[1]
{\noindent
{\em Proof of {#1}}.}
{\nopagebreak\mbox{}\hfill $\Box$\par\addvspace{0.5cm}}

\newcommand\mytop[2]{\genfrac{}{}{0pt}{}{#1}{#2}}

\newcommand{\eps}{\varepsilon}
\newcommand{\ov}{\overline}

\def\supp{\mathrm{supp}}

\def\id{\mathrm{id}}

\def\R{\mathbb{R}}

\def\dim{\mathrm{dim}}

\newcommand{\cB}{{\mathcal B}}
\newcommand{\cC}{{\mathcal C}}
\newcommand{\cD}{{\mathcal D}}

\newcommand{\cF}{{\mathcal F}}

\newcommand{\cH}{{\mathcal H}}

\newcommand{\cM}{{\mathcal M}}

\newcommand{\cR}{{\mathcal R}}

\newcommand{\al}{\alpha}

\newcommand{\de}{\delta}

\newcommand{\la}{\lambda}
\newcommand{\om}{\omega}
\newcommand{\si}{\sigma}
\newcommand{\Ga}{\Gamma}
\newcommand{\De}{\Delta}

\newcommand{\Om}{\Omega}

\newcommand{\pa}{\partial}
\def\id{\mathrm{id}}

\newcommand{\xbar}{\bar{x}}
\newcommand{\ybar}{\bar{y}}

\newcommand{\beq[1]}{\begin{equation}\label{eq:#1}}
\newcommand{\eeq}{\end{equation}}

\numberwithin{equation}{section}

\DeclareMathOperator{\dist}{dist}

\DeclareMathOperator{\Range}{Range}

\newcommand\set[1]{\left\{\,#1\,\right\}}  

\begin{document}

\title{The Morse property for functions of Kirchhoff-Routh path type}
\author{Thomas Bartsch\footnote{Supported by funds of ``Agreement between Sapienza University of Roma and University of Giessen''.} \and Anna Maria Micheletti \and Angela Pistoia}

\date{}
\maketitle


\begin{abstract}
For a bounded domain $\Omega\subset\mathbb{R}^n$ let $H_\Omega:\Omega\times\Omega\to\mathbb{R}$ be the regular part of the Dirichlet Green function for the Laplace operator. Given a fixed arbitrary ${\mathcal C}^2$ function $f:{\mathcal D}\to\mathbb{R}$, defined on an open subset ${\mathcal D}\subset\mathbb{R}^{nN}$, and fixed coefficients $\lambda_1,\dots,\lambda_N\in\mathbb{R}\setminus\{0\}$ we consider the function $f_\Omega:{\mathcal D}\cap\Om^N\to\R$ defined as
\[
  f_\Omega(x_1,\dots,x_N)
   = f(x_1,\dots,x_N) - \sum_{j,k=1}^N \lambda_j\lambda_k H_\Omega(x_j,x_k).
\]
We prove that $f_\Omega$ is a Morse function for {\em most domains} $\Omega$ of class ${\mathcal C}^{m+2,\al}$, any $m\ge0$, $0<\alpha<1$. This applies in particular to the Robin function $h:\Omega\to\mathbb{R}$, $h(x)=H_\Omega(x,x)$, and to the Kirchhoff-Routh path function where $\Omega\subset\mathbb{R}^2$, ${\mathcal D}=\{x\in\mathbb{R}^{2N}: \text{$x_j\ne x_k$ for $j\ne k$}\}$, and
\[
  f(x_1,\dots,x_N) = - \frac{1}{2\pi}\sum_{\genfrac{}{}{0pt}{}{j,k=1}{j\ne k}}^N\lambda_j\lambda_k\log|x_j-x_k|.
\]
\end{abstract}



\section{Introduction and main results}\label{sec:intro}
The paper is concerned with the Morse property of functions of the form
\begin{equation}\label{eq:f_Om}
  f_\Om(x_1,\dots,x_N) = f(x_1,\dots,x_N) - \sum_{j,k=1}^N \la_j\la_k H_\Om(x_j,x_k).
\end{equation}
Here $\Om\subset\R^n$ is a bounded domain, $H_\Om:\Om\times\Om\to\R$ is the regular part of the Green function for the Laplace operator with Dirichlet boundary conditions, and $f:\cD\to\R$ is a function of class $\cC^2$, defined on an open subset $\cD\subset\R^{nN}$. The function $f_\Om$ is then defined on $\cD\cap\Om^N$. Throughout the paper the function $f$ and the coefficients $\la_1,\dots,\la_N\in\R\setminus\{0\}$ are fixed arbitrarily. Our goal is to prove that for a generic domain $f_\Om$ is a Morse function, that is, all of its critical points are non-degenerate. We also have a symmetric version of our result.

Functions of the form \eqref{eq:f_Om} appear as singular limits in a variety of nonlinear partial differential equations. Most prominent is the Kirchhoff-Routh path function
\[
  \cH_{KR}(x_1,\dots,x_N) = - \frac{1}{2\pi}\sum_{\mytop{j,k=1}{j\ne k}}^N\la_j\la_k\log|x_j-x_k|
                               - \frac{1}{2\pi}\sum_{\genfrac{}{}{0pt}{}{j,k=1}{j\ne k}}^N\lambda_j\lambda_k\log|x_j-x_k|
\]
from fluid dynamics, introduced by Kirchhoff \cite{Kirchhoff:1876}, Routh \cite{Routh:1881} and Lin \cite{Lin:1941a, Lin:1941a}; see also \cite{Marchioro-Pulvirenti:1994, Newton:2001} for modern treatments. Here $\Om\subset\R^2$, $\cD=\{x=(x_1,\dots,x_N)\in\R^{2N}: \text{$x_j\ne x_k$ for $j\ne k$}\}$, and $\cH_{KR}=f_\Om$ with
\[
  f(x) = - \frac{1}{2\pi}\sum_{\mytop{j,k=1}{j\ne k}}^N\la_j\la_k\log|x_j-x_k|.
\]
The Hamiltonian system $\dot{x}_k=\nabla_{x_k}\cH_{KR}(x_1,\dots,x_N)$, $k=1,\dots,N$, describes the dynamics of $N$ point vortices with vortex strengths $\la_1,\dots,\la_N$ in an ideal fluid in $\Om$. Thus critical points of $\cH_{KR}$ are stationary point vortex solutions of the Euler equation in vorticity form. Knowing that these critical points are non-degenerate is very helpful for further investigations, for instance about the stability of the stationary solutions or the existence of periodic solutions near an equilibrium, or about the existence of heteroclinic or homoclinic solutions of the point vortex Hamiltonian system. It is also helpful for the desingularization of the point vortices, that is, for finding regular solutions of the Euler equation with vortices close to the singular point vortices.

Functions of the form \eqref{eq:f_Om} appear also as singular limits in a variety of nonlinear elliptic boundary value problems, for instance the renormalized energy for the Ginzburg-Landau equation. Other examples are Liouville type equations or mean field type equations. As in the case of the Euler equation, a non-degenerate critical point $(x_1,\dots,x_N)\in\cD\cap\Om^N$ of $f_\Om$ yields solutions of the elliptic problem that develop peaks (bubbles) at the points $x_1,\dots,x_N\in\Om$. Thus there is ample motivation for studying the Morse property of functions of the form $f_\Om$.

In order to formulate our result precisely, we fix a bounded domain $\Om\subset\R^n$ of class $\cC^{m+2,\al}$, $m\ge0$, $0<\al<1$. We also fix an open subset $\cD\subset\R^{nN}$, a $\cC^2$ function $f:\cD\to\R$ and parameters $\la_1,\dots,\la_N\in\R\setminus\{0\}$, and consider the function $f_\Om:\cD\cap\Om^N\to\R$ defined in \eqref{eq:f_Om}. The domain variations will be parameterized by elements from the Banach space $\cC^{m+2,\al}(\ov{\Om},\R^n)$ which is provided with the standard norm $\|\cdot\|_{m+2,\al}$.
For $\psi\in\cC^{m+2,\al}(\ov{\Om},\R^n)$ the set
\[
  \Om_\psi := (\id+\psi)(\Om) = \{x+\psi(x):x\in\Om\}
\]
is again a bounded domain of class $\cC^{m+2,\al}$ provided $\|\psi\|_{\cC^1}<\rho(\Om)$ is small. Setting
\[
  \cB^{m+2,\al}(\Om) := \set{\psi\in\cC^{m+2,\al}(\ov{\Om},\R^n): \|\psi\|_{\cC^1}<\rho(\Om)}
\]
we can now state our main result.

\begin{thm}\label{thm:main}
  Let $\Om\subset\R^n$ be a bounded domain of class $\cC^{m+2,\al}$ with $m\ge0$, $0<\al<1$. Then the set
  \[
    \cM^{m+2,\al}(\Om) := \set{\psi\in\cB^{m+2,\al}(\Om):\text{ $f_{\Om_\psi}$ is a Morse function}}
  \]
  is a residual (hence dense) subset of $\cB^{m+2,\al}(\Om)$.
\end{thm}

\begin{rem} \rm
a) Theorem~\ref{thm:main} applies in particular to the Kirchhoff-Routh path function $\cH_{KR}$. It also applies to the Robin function $h:\Om\to\R$, $h(x)=H_\Om(x,x)$. This case has already been treated in \cite{Micheletti-Pistoia:2014}. However, the proof in  \cite{Micheletti-Pistoia:2014} has a gap which is being fixed in this paper.

b) There are a number of results concerning the existence of critical points of the Kirchhoff-Routh path function; see \cite{Bartsch-DAprile-Pistoia:2013, Bartsch-Micheletti-Pistoia:2006, Bartsch-Pistoia:2015, Bartsch-Pistoia-Weth:2010, delPino-Kowalczyk-Musso:2005, Kuhl:2015, Kuhl:2016}, and references therein. In these and other papers the Kirchhoff-Routh path function appears as a singular limit when solving certain nonlinear elliptic boundary value problems. The non-degeneracy of the critical points is helpful when passing to the elliptic problem in that it often allows to replace a degree or variational argument by the contraction mapping principle, thus making the existence proof constructive. See also \cite{delPino-Kowalczyk-Musso:2005, Esposito-Grossi-Pistoia:2005} for applications to Liouville type equations and to mean field type equations.

c) Let us state some results on the dynamics of vortices where Theorem~\ref{thm:main} (and the symmetric version Theorem~\ref{thm:sym} below) are useful. In \cite{Lin-Lin:1997} the authors obtain solutions $u_\eps:\Om\to\R^2$ of the Ginzburg-Landau equation $-\De u=\frac{1}{\eps^2}(1-|u|^2)u$ in $\Om$ with Dirichlet boundary condition $u=g$ on $\pa\Om$ with vortices converging as $\eps\to0$ towards a prescribed non-degenerate critical point of the associated renormalized energy function, a function of the form $f_\Om$.

In \cite{Bartsch-Dai:2016, Bartsch-Gebhard:2016} periodic solutions of the $N$-vortex problem from fluid dynamics have been found near stable critical points of the Robin function. If the critical point is non-degenerate then it has been proved in \cite{Bartsch:2016} that there exists a smooth one-parameter curve of periodic solutions. The proof in \cite{Bartsch:2016} is based on the contraction mapping principle whereas the methods used in \cite{Bartsch-Dai:2016, Bartsch-Gebhard:2016} are non-constructive.

In \cite{Gebhard:2017}, for $N'>N$ non-stationary periodic solutions of the $N'$-vortex problem are constructed near non-degenerate critical points of the Kirchhoff-Routh path function (for $N$ vortices). Our result shows that the non-degeneracy assumption is generically true.
\end{rem}

If the domain $\Om$ is symmetric with respect to a subgroup $G\subset O(n)$ then $G$ also acts on $\R^{nN}$ via $g*(x_1,\dots,x_N)=(gx_1,\dots,gx_N)$. If moreover $\cD\subset\R^{nN}$ and $f:\cD\to\R$ are invariant under this action then also $f_\Om:\cD\cap\Om^N\to\R$ is invariant. In that case one can expect many critical points. For the Kirchhoff-Routh path function results in this direction can be found in \cite{Bartsch-Micheletti-Pistoia:2006, Bartsch-Pistoia-Weth:2010, Kuhl:2016}. There $\Om\subset\R^2$ is invariant under a finite group $G\subset O(2)$. A symmetric version of Theorem~\ref{thm:main} would therefore be useful where one only considers perturbations $\psi:\ov{\Om}\to\R^n$ from the set
\[
  \cB^{m+2,\al}_G(\Om) := \set{\psi\in\cB^{m+2,\al}(\Om): \text{ $\psi$ is equivariant}}.
\]
We can prove the following result.

\begin{thm}\label{thm:sym}
  Suppose $\Om$ is invariant under a finite subgroup $G\subset O(n)$ and $f:\cD\to\R$ is invariant with respect to the induced action of $G$ on $\cD\subset\R^n$. Then the set
  \[
    \cM^{m+2,\al}_G(\Om) := \set{\psi\in\cB^{m+2,\al}_G(\Om) :\text{ $f_{\Om_\psi}$ is a Morse function}}
  \]
  is a residual (hence dense) subset of $\cB^{m+2,\al}_G(\Om)$.
\end{thm}

As a corollary we obtain that the critical points of the Kirchhoff-Routh path function on symmetric domains found in \cite{Bartsch-Micheletti-Pistoia:2006, Bartsch-Pistoia-Weth:2010, Kuhl:2016} are nondegenerate for a generic symmetric domain.

It would be very interesting to allow symmetries with respect to compact subgroups $G\subset O(n)$. Observe that a critical point $x\in\cD\cap\Om^N$ of $f_\Om$ generates an orbit $Gx=\set{gx:g\in G}$ of critical points, hence critical points are always degenerate when $\dim(G)>0$. In that case one requires that all critical orbits $Gx$ are non-degenerate, i.e.\ the Hessian of $f_\Om$ is non-degenerate on the normal space to $Gx$. Unfortunately we could not deal with this case. The proof of Theorem~\ref{thm:main} is based on an abstract transversality theorem from \cite{Henry:2005}. Since we are not aware of an equivariant version of this theorem we only consider the case of finite groups $G$ here.

\section{Differentiability of $\mathbf H_\Om$ with respect to domain variations}

We fix a bounded $\cC^{m+2,\al}$ domain $\Om\subset\R^n$ with Green function $G_\Om=\Ga-H_\Om$. Here
\begin{equation}\label{Gamma}
 \Ga(x,y) = \begin{cases}
              -\frac{1}{2\pi}\ln|x-y| & \mbox{if $n=2$,} \\
              \frac{1}{(n-2)\om_n}|x-y|^{2-n} & \mbox{if $n\ge3$},
            \end{cases}
\end{equation}
is the singular part, and the regular part of the Green function. The regular part is the harmonic function with the same boundary values as the singular part, i.e.\ for any $y\in\Om$
\begin{equation}\label{H}
  \left\{\begin{aligned}
    &\De_x H_\Om(x,y)=0 \ &&\hbox{for}\ x\in\Om,\\
    &H_\Om(x,y)=\Ga(x,y) \ &&\hbox{for}\ x\in\pa\Om.\\
  \end{aligned}\right.
\end{equation}
In this section we show that $H_\Om$ is of class $\cC^1$ with respect to domain variations.

\begin{prop}\label{prop:hadamard}
  The map
  \[
    \cH_\Om:\Om\times\Om\times\cB^{m+2,\al}(\Om) \to \R, \quad \cH_\Om(x,y,\psi) = H_{\Om_\psi}(x+\psi(x), y+\psi(y))
  \]
  is of class $\cC^1$. Moreover, for $x,y\in\Om$ and $\phi\in\cC^k$ there holds:
  \begin{equation}\label{eq:hadamard}
    D_\psi\cH_\Om(x,y,0)[\phi] = \int_{\pa \Om}\langle \phi(z),\nu(z)\rangle\pa_{\nu_z}G_\Om(x,z)\pa_{\nu_z} G_\Om(y,z)\,d\si_z.
  \end{equation}
\end{prop}

In dimension $N=2$ the formula \eqref{eq:hadamard} goes back to Hadamard \cite{Hadamard:1908}.

\begin{proof}

It is clear that $\cH_\Om$ is $\cC^1$ in $(x,y)$. For $y\in\Om$ fixed we consider the map
\[
  \cH_{\Om,y}: \cB^{m+2,\al}(\Om)\to\cC^{m+2,\al}(\ov\Om,\R),\quad \cH_{\Om,y}(\psi)(x) = \cH_\Om(x,y,\psi) = H_{\Om_\psi}(x+\psi(x), y+\psi(y)).
\]

{\sc Step 1:} $\cH_{\Om,y}$ is continuous at $\psi=0$.\\
For $\phi\in\cB^{m+2,\al}(\Om)$ we write $\Phi=(\id +\phi)^{-1}$ and set $w_\phi(x):=H_{\Om_\phi}(x+\phi(x), y+\phi(y))$. Since $u:=w_\phi\circ\Phi$ is the unique solution of
\[
  \De u(z) =0 \quad\text{for $z\in\Om_\phi$},\qquad u(z)=\Ga(z,y+\phi(y))\quad\text{for $z\in\pa\Om_\phi$,}
\]
a straightforward computation shows that the map $w_\phi\in\cC^{m+2,\al}(\ov{\Om})$ is the unique solution of
\[
  \begin{cases}\displaystyle
    \sum_{i,j=1}^{n}a^{ij}_\phi(x)\frac{\pa^2w_\phi}{\pa x_i\pa x_j}(x)+\sum_{i=1}^{n}b^i_\phi(x)\frac{\pa w_\phi}{\pa x_i}(x) = 0 &\mbox{for }x\in\Om,\\
    w_\phi(x) = \Ga_\phi(x) & \mbox{for } x\in\pa\Om,
  \end{cases}
\]
with $a^{ij}_\phi(x)=\nabla\Phi_i\big(x+\phi(x)\big)\cdot\nabla\Phi_j\big(x+\phi(x)\big)$, $b^i_\phi(x)=\De\Phi_i\big(x+\phi(x)\big)$, and $\Ga_\phi(x)=\Ga\big(x+\phi(x),y+\phi(y)\big)$. It is not difficult to prove that $a^{ij}_\phi\to\de_{ij}$, $b^i_\phi\to0$ in $\cC^{m,\al}$, and $\Ga_\phi\to\Ga(\,\cdot\,,y)$ in $\cC^{m+2,\al}$, as $\phi\to0$ in $\cC^{m+2,\al}$. Standard elliptic estimates (see \cite[Theorem~6.6]{Gilbarg-Trudinger:1977}) imply that $w_\phi\to H_\Om(\,\cdot\,,y)$ in $\cC^{m+2,\al}$.

{\sc Step 2:} $\cH_{\Om,y}$ is continuous.\\
Observe that
\begin{equation}\label{eq:step2}
  \cH_{\Om,y}(\psi+\phi)(x) = \cH_{\Om_\psi,y+\psi(y)}\left(\phi\circ(\id+\psi)^{-1}\right)\big(x+\psi(x)\big).
\end{equation}
Applying {\sc Step 1} with $\Om_\psi$ instead of $\Om$ and $y+\psi(y)$ instead of $y$ we obtain that the map $\cH_{\Om_\psi,y+\psi(y)}$ is continuous at $0$, hence $\cH_{\Om,y}$ is continuous at $\psi$.

{\sc Step 3:} $\cH_{\Om,y}$ is Gateaux differentiable.\\
Using \eqref{eq:step2} and arguing as in {\sc Step 2} it is sufficient to show that $\cH_{\Om,y}$ is Gateaux differentiable at $\psi=0$. A straightforward computation gives for $\phi\in\cB^{m+2,\al}(\Om)$:
\[
  \begin{aligned}
    \De_x\cH_{\Om,y}(\phi)(x)
      & = 2\sum^n_{i,j=1}\frac{\pa^2}{\pa z_i\pa z_j}\Big| _{z=x+\phi(x)} H_{\Om_\phi}(z,y+\phi(y))\frac{\pa}{\pa x_i}\phi_j(x)\\
      &\hspace{1cm} + \sum^n_{i,j=1}\frac{\pa^2}{\pa z_i\pa z_j}\Big|_{z=x+\phi(x)} H_{\Om_\phi}(z,y+\phi(y)) \nabla\phi_i(x)\cdot\nabla\phi_j(x)\\
      &\hspace{1cm} +\nabla_z\big| _{z=x+\phi(x)}H_{\Om_\phi}(z,y+\phi(y))\cdot \De\phi(x)
 \end{aligned}
\]
It follows that $\displaystyle w=\lim\limits_{\eps\to 0}\frac{1}{\eps}\big(\cH_{\Om,y}(\eps\phi)-\cH_{\Om,y}(0)\big)$ satisfies
\[
  \De w(x) = 2\sum_{i,j=1}^N\frac{\pa^2}{\pa x_i\pa x_j}H_\Om(x,y)\frac{\pa\phi_j}{\pa x_i} (x) + \nabla_x H_\Om(x,y)\cdot \De \phi (x)
    \quad\text{for $x\in\Om$.}
\]
Moreover, one easily sees that
\[
  w(x) = -\frac{(x-y)\cdot(\phi(x)-\phi(y))}{\omega_n|x-y|^n} \quad\text{for $x\in\pa\Om$.}
\]
This implies:
\[
  \begin{aligned}
  w(x)
   &= -\int_\Om \left(2\sum_{i,j=1}^n\frac{\pa^2}{z_i\pa z_j}H_\Om(z,y)\frac{\pa\phi_j}{\pa z_i}(z)+\nabla_z H_\Om(z,y)\cdot\De\phi(z)\right)G(x,z)\,dz\\
   &\hspace{1cm}
        + \frac {1}{\om_n}\int_{\pa\Om}\frac{(z-y) \cdot \big(\phi(z)-\phi (y)\big)}{|z-y|^n} \pa_{\nu_z}G(x,z)\,d\si_z
  \end{aligned}
\]
Now we compute, using $G(x,z)=0$ for $z\in\pa\Om$ and $\De_z H_\Om(z,y)=0$:
\[
  \begin{split}
  \begin{aligned}
    &\sum_{i,j=1}^n\int_\Om\frac{\pa^2}{\pa z_i \pa z_j} H_\Om(z,y) \frac{\pa\phi_j}{\pa z_i}(z) G(x,z)\,dz\\
    &\hspace{.5cm}
     = \sum_{i,j=1}^n\int_\Om\frac{\pa}{\pa z_i}\left(\frac{\pa ^2}{\pa z_i\pa z_j} H_\Om(z,y) \phi_j(z) G(x,z)\right)\,dz\\
    &\hspace{1cm}
         -\sum_{i,j=1}^n\int_\Om \frac{\pa}{\pa z_i}\left(\frac{\pa^2}{\pa z_i \pa z_j} H_\Om(z,y) G(x,z)\right)\phi_j(z)\,dz\\
    &\hspace{.5cm}
     = \sum_{i,j=1}^n\int_{\pa\Om}\frac{\pa ^2}{\pa z_i \pa z_j}H_\Om(z,y) \phi_j(z) G(x,z) \nu_i (z)\,d\si_z\\
    &\hspace{1cm}
         -\sum_{j=1}^n\int_\Om \frac{\pa}{\pa z_j}\De_z H_\Om (z,y)  G(x,z)\phi_j(z)\,dz
         - \sum_{i,j=1}^n\int_\Om\frac{\pa ^2}{\pa z_i \pa z_j}H_\Om(z,y)\frac{\pa}{\pa z_i} G(x,z) \phi_j (z)\,dz\\
    &\hspace{.5cm}
     = - \sum_{i,j=1}^n\int_\Om\frac{\pa ^2}{\pa z_i \pa z_j}H_\Om(z,y)\frac{\pa}{\pa z_i} G(x,z) \phi_j (z)\,dz
  \end{aligned}
  \end{split}
\]
Similarly we obtain:
\[
  \begin{split}
  \begin{aligned}
     &\int_\Om \nabla_z H(z,y) \cdot \De\phi(z) G(x,z) dz
      = \sum_{j=1}^n\int_\Om\frac{\pa}{\pa z_j} H_\Om(z,y)\De\phi_j(z)G(x,z) dz\\
     &\hspace{.5cm}
      =\sum_{i,j=1}^n\int_\Om\frac{\pa}{\pa z_i} \left(\frac{\pa}{\pa z_j} H_\Om(z,y)\frac{\pa}{\pa z_i}\phi_j(z) G(x,z)\right) dz\\
     & \hspace{1cm}
        -\sum_{i,j=1}^n\int_\Om\frac{\pa}{\pa z_i}\left(\frac{\pa}{\pa z_j} H_\Om(z,y) G(x,z)\right)\frac{\pa}{\pa z_i}\phi_j (z) dz\\
     &\hspace{.5cm}
      = -\sum_{i,j=1}^n\int_\Om\frac{\pa}{\pa z_i}\left(\frac{\pa}{\pa z_i}\left(\frac{\pa}{\pa z_j}H_\Om(z,y) G(x,z)\right)\phi_j(z)\right)dz\\
     &\hspace{1cm}
         +\sum_{i,j=1}^n\int_\Om\frac{\pa^2}{\pa z^2_i}\left(\frac{\pa}{\pa z_j}H_\Om (z,y)G(x,z)\right)\phi_j(z) dz\\
     &\hspace{.5cm}
      = -\sum_{i,j=1}^n\int_{\pa\Om}\frac{\pa}{\pa z_i}\left(\frac{\pa}{\pa z_j}H_\Om(z,y) G(x,z)\right)\phi_j(z)\nu_i(z) d\si_z\\
     &\hspace{1cm}
         +2\sum_{i,j=1}^n\int_\Om\frac{\pa^2}{\pa z_i\pa z_j}H_\Om(z,y)\frac{\pa}{\pa z_i}G(x,z)\phi_j(z) dz\\
  \end{aligned}
  \end{split}
\]
If $\phi(y)=0$ it follows that
\[
  \begin{split}
  \begin{aligned}
    w(x)
     &= \sum_{i,j=1}^n\int_{\pa\Om} \frac{\pa}{\pa z_i} \left(\frac{\pa}{\pa z_j} H_\Om (z,y)G(x,z)\right)\phi_j(z) \nu_i(z)\,d\si_z\\
     &\hspace{1cm}
          +\frac{1}{\om_n}\int_{\pa\Om} \frac{(z-y)\cdot\big(\phi(z)-\phi(y)\big)}{|z-y|^n} \pa_{\nu_z}G(x,z)\,d\si_z\\
     &= \sum_{i,j=1}^n\int_{\pa\Om} \frac{\pa}{\pa z_j}H_\Om(z,y)\frac{\pa}{\pa z_i} G(x,z)\phi_j(z) \nu_i(z)\,d \si_z\\
     &\hspace{1cm}
          +\frac{1}{\om_n}\int_{\pa\Om} \frac{(z-y)\cdot\big(\phi(z)-\phi(y)\big)}{|z-y|^n} \pa_{\nu_z}G(x,z)\,d\si_z\\
     &= \int_{\pa\Om}\nabla_z H_\Om(z,y)\cdot \phi(z)\pa_{\nu_z} G(x,z)\,d\si_z\\
     &\hspace{1cm}
          +\frac{1}{\om_n}\int_{\pa\Om}\frac{(z-y)\cdot(\phi(z)-\phi(y))}{|z-y|^n} \pa_{\nu_z}G(x,z)\,d\si_z\\
     &= -\int_{\pa\Om}\nabla_z G(z,y)\cdot \phi(z)\cdot \pa_{\nu_z} G(x,z)\,d\si_z\\
     &= -\int_{\pa\Om}\langle \nu(z), \phi(z)\rangle \pa_{\nu_z}G(z,y)\pa_{\nu_z}G(x,z)\,d\si_z
  \end{aligned}
  \end{split}
\]
Here we used $\nabla_z G(z,y)=\pa_{\nu_z}G(z,y)\cdot \nu(z)$ for $z\in\pa\Om$. Thus we have proved the Gateaux differentiability at $\psi=0$ in the direction $\phi$, and equation~\eqref{eq:hadamard}, provided $\phi(y)=0$. Since $\cH(x,y,\phi)$ and \eqref{eq:hadamard} are symmetric in $x$ and $y$, equation~\eqref{eq:hadamard} also holds if $\phi(x)=0$. Now a general $\phi\in\cC^{m+2,\al}(\ov\Om,\R^n)$ can be written as $\phi=\phi_1+\phi_2$ with $\phi_1,\phi_2\in\cC^{m+2,\al}(\ov\Om,\R^n)$ and such that $\phi_1(x)=0$ and $\phi_2(y)=0$. Therefore $\cH_{\Om,y}$ is Gateaux differentiable at $\psi=0$ in any direction $\phi\in\cC^{m+2,\al}(\ov\Om,\R^n)$.

{\sc Step 4:} $\cH_{\Om,y}$ is continuously Frechet differentiable.\\
Using \eqref{eq:hadamard} and \eqref{eq:step2} we deduce for the Gateaux derivative at $\psi$ in the direction $\phi$:
\[
  \begin{aligned}
  D_\psi\cH_{\Om,y}(\psi)[\phi](x)
   &= D\cH_{\Om_\psi,y+\psi(y)}(0)[\phi\circ(\id+\psi)^{-1}](x+\psi(x)) \\
   &= \int_{\pa\Om_\psi}\langle \phi\circ(\id+\psi)^{-1}(z),\nu(z)\rangle\pa_{\nu_z}G_{\Om_\psi}(x+\psi(x),z)\pa_{\nu_z} G_{\Om_\psi}(y+\psi(y),z)\,d\si_z.
  \end{aligned}
\]
Making the transformation $\zeta=(\id+\psi)^{-1}(z)$ and using {\sc Step 1} one sees that the Gateaux derivative of $\cH_{\Om,y}$ is continuous in $x$ and $\psi$.

{\sc Step 5:} $\cH_\Om$ is continuously Frechet differentiable.\\
By {\sc Step 4} $\cH_\Om(x,y,\psi)$ is continuously Frechet differentiable in $x$ and $\psi$. The claim follows easily using the symmetry  $\cH_\Om(x,y,\psi)=\cH_\Om(y,x,\psi)$.
\end{proof}

As a corollary we obtain the differentiability of the Robin function with respect to domain perturbations.

\begin{cor}
  The map
  \[
    \cR_\Om:\Om\times\cB^{m+2,\al}(\Om) \to \R, \quad \cR_\Om(x,\psi) = H_{\Om_\psi}(x+\psi(x), x+\psi(x))
  \]
  is of class $\cC^1$. Moreover, for $x\in\Om$ and $\phi\in\cC^{m+2,\al}(\ov\Om,\R^n)$ there holds:
  \begin{equation}\label{eq:hadamard_robin}
    D_\psi\cR_\Om(x,0)[\phi] = 2\int_{\pa \Om}\big\langle \phi(z),\nu(z)\big\rangle|\pa_{\nu_z}G(x,z)|^2\,d\si_z.
  \end{equation}
\end{cor}

\section{Proof of Theorems \ref{thm:main} and \ref{thm:sym}}
The proof is based on the following theorem which is a special case of \cite[Theorem~5.4]{Henry:2005}.

\begin{thm}\label{thm:trans}
 Let $X,Y,Z$ be three Banach spaces and let $\cF:U\to Z$ be a $C^1$ map defined on an open subset $U\subset X\times Y$. Assume that:
 \begin{itemize}
 \item[(i)] for any $(\xbar,\ybar)\in \cF^{-1}(0)$, the map $\frac{\pa \cF}{\pa x}(\xbar,\ybar):X\to Z$ is a Fredholm operator of index $0$;
 \item[(ii)] $0$ is a regular value of $\cF$, i.e.\ the operator $D\cF(\xbar,\ybar):X\times Y\to Z$ is onto at every point $(\xbar,\ybar)\in\cF^{-1}(0)$;
 \item[(iii)] the map  $\pi\circ i:\cF^{-1}(0)\subset X\times Y\to Y$ is $\si$-proper, i.e.\ $\cF^{-1}(0)=\bigcup_{j=1}^{+\infty} M_j$ is a countable union of sets $M_j$ and the restriction $\pi\circ i|_{M_j}$ is proper for any $j$. Here $i:\cF^{-1}(0)\to X\times Y$ is the inclusion and $\pi:X\times Y\to Y$ the projection.
 \end{itemize}
 Then the set
 \[
   Y_{reg}:=\set{y\in Y\ :\ 0\ \hbox{is  a regular value of } \cF(\cdot,y)}
 \]
 is a  residual subset of $Y$, i.e.\ $Y\setminus Y_{reg}$ is a countable union of closed subsets without interior points.
\end{thm}

Observe that $\cF(\cdot,y)$ is defined on the set $U_y=\{x\in X:(x,y)\in U\}$. If $U_y=\emptyset$ then $y\in Y_{reg}$.

\begin{altproof}{Theorem~\ref{thm:main}}
We apply Theorem~\ref{thm:trans} in the following setting.  Let $X=Z=\R^{nN}$, $Y=\cC^{m+2,\al}(\ov{\Om},\R^n)$ and set
\[
  U := \set{(x,\psi)\in\R^{nN}\times\cB^{m+2,\al}(\Om): \big(x_1+\psi(x_1),\dots,x_N+\psi(x_N)\big)\in \cD}.
\]
Consider the map $\cF_\Om:U\to\R^{nN}$ defined by
\[
  \begin{aligned}
  \cF_\Om(x,\psi)
   &= \nabla_x f_{\Om_\psi}(x_1,\dots,x_N)\\
   &= \nabla_x\left(f\big(x_1+\psi(x_1),\dots,x_N+\psi(x_N)\big) - \sum_{j,k=1}^N \la_j\la_k H_{\Om_\psi}\big(x_j+\psi(x_j),x_k+\psi(x_k)\big)\right).
  \end{aligned}
\]
For the proof of Theorem~\ref{thm:main} we have to show that $\cM^{m+2,\al}(\Om)=Y_{reg}$ is residual in $Y$.

{\sc Step 1:} $\cF_\Om$ satisfies (i) and (iii) from Theorem~\ref{thm:trans}.\\
Since $\dim X=\dim Z<\infty$ property (i) is trivially satisfied. In order to prove (iii) we set $\Om_j:=\set{x\in\Om:\dist(x,\pa\Om)\ge 1/j}$ and
\[
  U_j:=\set{(x,\psi)\in U:x_k\in\Om_j\text{ for $k=1,\dots,N$},\ \|\psi\|_{m+2,\al}\le \rho(\Om)-1/j}.
\]
Then the restriction $\pi\circ i|_{M_j}$ of $\pi\circ i$ to $M_j:=U_j\cap\cF_\Om^{-1}(0)$ is proper because $\Om_j$ is compact and $U_j$ is closed in $\R^{nN}\times\cB^{m+2,\al}(\Om)$. Clearly we have $\cF_\Om^{-1}(0)=\bigcup_{j=1}^\infty M_j$.

{\sc Step 2:} For every $\xbar\in\cD\cap\Om^N$ the operator $D_\psi\cF_\Om(\xbar,0):\cC^k(\ov{\Om},\R^n) \to \R^{nN}$ is onto.\\
Given $\xbar\in\cD\cap\Om^N$ we compute
\[
  \frac{\pa}{\pa\eps}\Big|_{\eps=0} \Big(F_{\Om_{\eps\phi}}\big(x_1+\eps\phi(x_1),\dots,x_N+\eps\phi(x_N)\big)\Big)
\]
for $\phi\in\cC^{m+2,\al}(\ov{\Om},\R^n)$ with $\xbar_1,\dots,\xbar_N\notin\supp(\phi)$. This last condition implies
\[
  \frac{\pa}{\pa\eps}\Big|_{\eps=0} \Big(f\big(x_1+\eps\phi(x_1),\dots,x_N+\eps\phi(x_N)\big)\Big) = 0
\]
for $x$ near $\xbar$. Therefore we have for $x$ near $\xbar$:
\[
  \begin{aligned}
   &\frac{\pa}{\pa\eps}\Big|_{\eps=0} \Big(F_{\Om_{\eps\phi}}\big(x_1+\eps\phi(x_1),\dots,x_N+\eps\phi(x_N)\big)\Big)\\
   &\hspace{1cm}
      = -\sum_{j,k=1}^N \la_j\la_k \frac{\pa}{\pa\eps}\Big|_{\eps=0} \left(H_{\Om_{\eps\phi}}(x_j+\eps\phi(x_j),x_k+\eps\phi(x_k))\right)\\
   &\hspace{1cm}
      = -\sum_{j,k=1}^N \la_j\la_k \int_{\pa \Om}\langle \phi(z),\nu(z)\rangle\pa_{\nu_z}G_\Om(\xbar_j,z)\pa_{\nu_z} G_\Om(\xbar_k,z)\,d\si_z.
  \end{aligned}
\]
When passing to the gradient
\[
  D_\psi\cF_\Om(\xbar,0)[\phi]
   = \nabla_x\Big|_{x=\xbar} \frac{\pa}{\pa\eps}\Big|_{\eps=0} \Big(F_{\Om_{\eps\phi}}\big(x_1+\eps\phi(x_1),\dots,x_N+\eps\phi(x_N)\big)\Big)
\]
it is useful to identify $\R^{nN}$ with $\R^N\otimes\R^n$. An element $(x_1,\dots,x_N)\in(\R^n)^N=\R^{nN}$ corresponds to $\sum_{k=1}^{N}e_k\otimes x_k$; here $e_1,\dots,e_N\in\R^N$ is the standard basis. With this notation we have:
\[
  \begin{aligned}
  D_\psi\cF_\Om(\xbar,0)[\phi]
   &= \nabla_x\Big|_{x=\xbar}\frac{\pa}{\pa\eps}\Big|_{\eps=0} \Big(F_{\Om_{\eps\phi}}\big(x_1+\eps\phi(x_1),\dots,x_N+\eps\phi(x_N)\big)\Big)\\
   &= -2\sum_{j,k=1}^N \la_j\la_k \int_{\pa \Om}\langle \phi(z),\nu(z)\rangle\pa_{\nu_z}G_\Om(\xbar_j,z)
       \big(e_k\otimes\nabla_{x_k}\pa_{\nu_z}G_\Om(\xbar_k,z)\big)\,d\si_z.
  \end{aligned}
\]
Let $v\in\big(\Range D_\psi\cF_\Om(\xbar,0)\big)^\perp\subset\R^{nN}$ be an arbitrary element of the orthogonal complement of the range of $D_\psi\cF_\Om(\xbar,0)$. We shall show that $v=\sum_{k=1}^{N}e_k\otimes v_k=0$, thus proving the claim. For every $\phi\in\cC_k$ with $\xbar_1,\dots,\xbar_N\notin\supp(\phi)$ there holds:
\begin{equation}\label{eq:v=0-1}
  \begin{aligned}
    0 &= \big\langle D_\psi\cF_\Om(\xbar,0)[\phi],v\big\rangle\\
      &= -2\sum_{j,k=1}^N \la_j\la_k \int_{\pa \Om}\langle \phi(z),\nu(z)\rangle\pa_{\nu_z}G_\Om(\xbar_j,z)\big\langle\nabla_{x_k}\pa_{\nu_z}G_\Om(\xbar_k,z),v_k\big\rangle\,d\si_z.
  \end{aligned}
\end{equation}
Since $\phi$ can be arbitrary on the boundary $\pa\Om$ we deduce for every $z\in\pa\Om$:
\begin{equation}\label{eq:v=0-2}
  \begin{aligned}
  0 &= \sum_{j,k=1}^N \la_j\la_k \pa_{\nu_z}G_\Om(\xbar_j,z)\big\langle\nabla_{x_k}\pa_{\nu_z}G_\Om(\xbar_k,z),v_k\big\rangle\\
    &= \left(\sum_{j=1}^N \la_j\pa_{\nu_z}G_\Om(\xbar_j,z)\right)
       \left(\sum_{k=1}^N \la_k\big\langle\nabla_{x_k}\pa_{\nu_z}G_\Om(\xbar_k,z),v_k\big\rangle\right)\\
    &= \pa_{\nu_z}h_1(z)\pa_{\nu_z}h_2(z)
  \end{aligned}
\end{equation}
where the functions $h_1,h_2:\ov{\Om}\setminus\{\xbar_1,\dots,\xbar_N\} \to \R$ are defined by
\[
  h_1(z) = \sum_{j=1}^N \la_j G_\Om(\xbar_j,z)
  \qquad\text{and}\qquad
  h_2(z) = \sum_{k=1}^N \la_k \big\langle\nabla_{x_k}G_\Om(\xbar_k,z),v_k\big\rangle\,.
\]
Observe that both $h_1$ and $h_2$ are harmonic for $z\in\Om$ and identically $0$ for $z\in\pa\Om$. Since $h_1$ is not identically equal to $0$ the unique continuation principle implies that the set $\{z\in\pa\Om:\pa_{\nu_z}h_1(z)=0\}$ does not contain an open subset of $\pa\Om$. Now \eqref{eq:v=0-2} implies that $\pa_{\nu_z}h_2(z)=0$ for all $z\in\pa\Om$. Using the unique continuation principle once more we deduce that $h_2\equiv0$. This implies $v=0$ because if $v_k\ne0$ then
\[
  h_2(\xbar_k+tv_k) = - \frac{\la_k}{\om_n t^{n-1}|v_k|^{n-2}} + O(1) \qquad\text{as $t\to0$.}
\]

{\sc Step 3:} $\cF_\Om$ satisfies (ii) from Theorem~\ref{thm:trans}.\\
This follows as in the proof of Proposition~\ref{prop:hadamard}. Simply observe for $(\xbar,\bar{\psi})\in U$ that
\[
  D_\psi\cF_\Om\big((\xbar_1,\dots,\xbar_N),\bar{\psi}\big)[\phi]
   = D_\psi\cF_{\Om_{\bar{\psi}}}\big(\xbar_1+\bar{\psi}(\xbar_1),\dots,\xbar_N+\bar{\psi}(\xbar_N),0\big)
       [\phi\circ(\id+\bar{\psi})^{-1}].
\]
Now apply {\sc Step 2} to $D_\psi\cF_{\Om_{\bar{\psi}}}(\xbar_1+\bar{\psi}(\xbar_1),\dots,\xbar_N+\bar{\psi}(\xbar_N),0)$ instead of $D_\psi\cF_\Om(\xbar_1,\dots,\xbar_N,0)$.
\end{altproof}

\begin{altproof}{Theorem~\ref{thm:sym}}
The proof extends to the equivariant setting in a straightforward way. The main point is to check that for \eqref{eq:v=0-2} to be true it is sufficient that \eqref{eq:v=0-1} holds for all $\phi\in\cC^k_G$ with $\xbar_1,\dots,\xbar_N\notin\supp(\phi)$. This requires the construction of equivariant test functions $\phi$ which we leave to the reader.
\end{altproof}

{\bf Acknowledgement.}  T.~B.\ thanks Universit\`a di Roma ``La Sapienza", in particular the Dipartimento SBAI, for its hospitality.

\vspace{2mm}\noindent
{\sc Thomas Bartsch}\\
 Mathematisches Institut\\
 Universit\"at Giessen\\
 Arndtstr.\ 2\\
 35392 Giessen, Germany\\
 Thomas.Bartsch@math.uni-giessen.de

\vspace{2mm}\noindent
{\sc Anna Maria Micheletti}\\
 Dipartimento di Matematica\\
 Universit\`a di Pisa\\
 Via Bonanno 25B\\
 56126 Pisa, Italy\\
 a.micheletti@dma.unipi.it

\vspace{2mm}\noindent
{\sc Angela Pistoia}\\
 Dipartimento SBAI\\
 Universit\`a di Roma ``La Sapienza"\\
 via Antonio Scarpa 16\\
 00161 Roma, Italy\\
 angela.pistoia@uniroma1.it

\end{document}